\documentclass[a4paper,11pt]{article}
\usepackage{amsfonts,amsmath,amssymb,graphicx,float}

\def\kaxxa{{\vcenter {\hrule height.2mm
\hbox{\vrule width .2mm height 2mm \kern 2mm
\vrule width .2mm} \hrule height .2mm}}}

\def\lesta{\hfill $\kaxxa$ \medskip}

\newcommand{\ern}{\mbox{ern}}

\newtheorem{Theorem}{Theorem}[section]
\newtheorem{Conjecture}{Conjecture}[section]
\newtheorem{Lemma}{Lemma}[section]
\newtheorem{Proposition}{Proposition}[section]
\newtheorem{Result}{Result}[section]
\newtheorem{Observation}{Observation}[section]
\newtheorem{Definition}{Definition}

\newtheorem{Example}{Example}[section]
\newtheorem{Corollary}{Corollary}[section]

\title{The degree-associated edge-reconstruction number of disconnected graphs and trees}
\author{Kevin J. Asciak\footnote{kevin.j.asciak@um.edu.mt} \\ Department of Mathematics \\ University of Malta \\ Malta}
\begin{document}

\maketitle

\begin{abstract}
An \emph{edge-card} of a graph $G$ is a subgraph formed by deleting an edge. The \emph{edge-reconstruction number} of a graph $G$, $ern(G)$, is the minimum number of edge-cards required to determine $G$ up to isomorphism. A \emph{da-ecard} is an edge-card which also specifies the degree of the deleted edge, that is, the number of edges adjacent to it. The \emph{degree-associated edge-reconstruction number}, $dern(G)$ is the minimum number of da-ecards that suffice to determine the graph $G$.
In this paper we state some known results on the edge-reconstruction number of disconnected graphs and trees. Then we investigate how the degree-associated edge-reconstruction number of disconnected graphs and trees vary from their respective edge-reconstruction number. We show how we can select two da-ecards to identify caterpillars uniquely. We also show that while $dern(tP_n) = 2 $ for $ n > 3$,  $dern(tP_3) = 3$ where $P_n$ is the path on $n$ vertices, and that, although $dern(K_{1,n}) = 1$, $dern(S_{p+1}^n) = 2 $ where $S_{p+1}^n$ is a tree obtained from the star $K_{1,n}$  by subdividing each edge $p$ times. Finally we conjecture that for any tree $T$, $dern(T) \leq 2 $.
\end{abstract}

Keywords: Degree-associated edge-reconstruction number, trees, disconnected graphs.

AMS Subject classification: 05C60

\section{Introduction}
All graphs are assumed to be simple, finite and undirected, and any graph-theoretic notations and definitions not explicitly defined can be found in \cite {bm} or \cite {ls}.

A \emph{vertex-deleted subgraph} $G - v$ is the unlabelled graph obtained by deleting, from the graph $G$, a vertex $v$ and all edges incident to $v$. The \emph{deck} of $G$, denoted $\mathcal{D}(G)$, is the multiset of  vertex-deleted subgraphs of $G$ and each member of $\mathcal{D}(G)$ is referred to as a \emph{card}. Our main focus in this paper will be on the analogously defined \emph{edge-cards} of $G$ which are edge-deleted subgraphs $G -e$ of $G$. The collection of the edge-cards is called the \emph{edge-deck} of $G$, denoted by $\cal{ED}(G)$.

The Reconstruction Conjecture proposed in 1942 by Kelly \cite{Ke42} and Ulam \cite{Ul60} states that a graph with at least three vertices, is uniquely determined, up to isomorphism, from its collection of vertex-deleted subgraphs. The most natural variation of the Reconstruction Conjecture (for a recent survey see \cite {l2}) is an analogue for deletion of edges. This is the Edge-Reconstruction Conjecture [Harary, 1964] which states that all graphs on at least four edges are edge reconstructible \cite{Har64}.

A \emph{reconstruction} (\emph{edge-reconstruction})\index{reconstruction  of a graph}\index{edge-reconstruction of a graph} of $G$ is a graph $H$ with $\cal{D}(G) = \cal{D}(H)$ ($\cal{ED}(G) = \cal{ED}(H)$).  A graph $G$ is \emph{reconstructible} (\emph{edge-reconstructible})\index{reconstructible   graph}\index{edge-reconstructible  graph} if every reconstruction of $G$ is isomorphic to $G$.  This means that $G$ is reconstructible (edge-reconstructible) if it can be obtained uniquely, up to isomorphism, from its deck (edge-deck).

For a reconstructible graph $G$, Harary and Plantholt \cite {hpl2} introduced the notion of the \emph{reconstruction number} of $G$, denoted by $rn(G)$, which is defined as the least number of vertex-deleted subgraphs of $G$ required in order to identify G uniquely; that is, $rn (G)$ is the size of the smallest subcollection of the deck of $G$ which is not contained in any other deck of another graph $H$ where $H \not\simeq  G$. It can be considered as a measure of the level of difficulty in reconstructing $G$ uniquely.
The simplest observation we can make is that $rn (G) \geq 3$. Reconstruction numbers are now known for various classes of graphs such as disconnected and regular graphs and trees \cite {aflm}.

A variation of reconstruction numbers is the \emph{class-reconstruction numbers}. Let $C$ be a class of graphs and let $G \in C$. Then the class-reconstruction number $Crn (G)$ is defined as the least number of vertex-deleted subgraphs required to determine $G$ from any other graph in $C$. Non-trivial class reconstruction numbers which have been studied include maximal planar graphs, trees and unicyclic graphs \cite {lh2,lh1,hpl1}.


Reconstructing the graph from the deck seems to be more difficult than reconstructing the graph from its edge-deck since more of the graph is left in an edge-deleted subgraph than in a vertex-deleted subgraph. However, it sometimes happens that more edge-deleted subgraphs are required for unique reconstruction than vertex-deleted subgraphs.

Motivated by the falsity of the Reconstruction Conjecture for directed graphs, Ramachandran \cite {r1,r2} weakened the Reconstruction Conjecture by considering the degree of the deleted vertex along with each vertex-deleted subgraph in a\emph{ degree-associated card}, referred to as \emph{da-card}.

The \emph{degree} of a vertex $v$ is the number of edges of $G$ incident to $v$.  A vertex of degree 0 is called an \emph{isolated vertex} and a vertex of degree 1 is called an \emph{endvertex}. The \emph{minimum degree} of a graph  $G$, denoted by $\delta(G)$, is the smallest number of edges incident to any vertex $v$ in $G$.

A da-card denoted by $( G - v, d)$  consists of a card $G - v$ in the deck of $G$ and the degree $d$ in $G$ of the deleted vertex $v$. The \emph{da-deck} is the multiset of da-cards. Ramachandran defined the \emph{degree-associated reconstruction number} of a graph $G$, denoted by $drn (G)$, to be the minimum number of da-cards necessary to determine $G$ uniquely.  Clearly $drn (G)$ is equivalent to the class reconstruction number of $G$ given that $G$ is in $C_m$, the class of graphs on $m$ edges.

The \emph{edge-reconstruction number}, the \emph{class edge-reconstruction number} and the \emph{degree-associated edge-reconstruction number} are analogously defined. 
A  \emph{degree-associated edge-card} or \emph{da-ecard} is a pair $( G - e, d(e))$ consisting of an edge-card $G - e $ in the edge deck of the graph $G$ and the degree of the edge $e$, denoted by $d(e)$, which is the number of edges adjacent to $e$, that is, $d(u) + d(v) - 2$ where $e = uv$. The multiset of all da-ecards is the \emph{da-edeck}.

Monikandan and Raj \cite {mdr,mr} initiated the study of the degree-associated edge-reconstruction number, denoted by $dern (G)$ which is the minimum $k$ such that some multiset of $k$ da-ecards determines $G$. Clearly $dern (G) \leq ern (G)$.
They determined $dern (G)$ where $G$ is a regular graph, a complete bipartite graph, a path, a wheel or a double star.  They also proved that $dern(G) \leq 2$ where $G$ is a complete 3-partite graph whose part-sizes differ by at most 1.  They showed that if $G$ is a graph obtained from $K_{1,m}$ by subdividing each edge at most once, then $dern(G) \leq 2$.

In her study, Myrvold \cite{my} proposed the \emph{adversary reconstruction number} of $G$, denoted by $adv\text{-} rn(G)$, which is the smallest value of $k$ such that no subdeck of $G$ containing $k$ cards is in the deck of any other graph which is not isomorphic to $G$. Therefore $adv\text{-} rn(G)$ equals 1 plus the largest number of cards which $G$ has in common with any graph not isomorphic to it. The \emph{adversary edge-reconstruction number} of $G$, denoted by $adv\text{-} ern(G)$, is analogously defined. Recently Monikandan et al. \cite{mrjs} also introduced a similar parameter called the \emph{adversary degree-associated edge-reconstruction number} of a graph $G$, denoted by $adv\text{-} dern(G)$, which is the least number $k$ such that every collection of $k$ da-ecards of $G$ is not contained in the da-edeck of any other graph $H$ such that $H \not\simeq G $. From the definitions, it follows that $ern(G) \leq adv\text{-} ern(G)$ and $dern(G) \leq min \{ ern(G), adv\text{-} dern(G) \}$. Moreover, if all the da-ecards of a graph $G$ are isomorphic, then $dern(G) = adv\text{-} dern(G)$.

Ma, Shi, Spinoza and West \cite {mssw} recently showed also that for all complete multipartite graphs and their complements $dern$ is usually 2 except for some exceptions. They also pointed out that a significant difference between vertex and edge degree-associated reconstruction number is that while trivially, a graph and its complement have the same $drn$ \cite {bw}, they need not have the same value of $dern$.

In this paper we first state some known results on the edge-reconstruction number of disconnected graphs and trees and then present some new results on their respective degree-associated edge-reconstruction numbers. There is a large gap between the value of $ern(G) =  3$ for disconnected graph $G$ with at least two non-isomorphic components and the value $ern(G) = t + 2$ for disconnected graph $G$ with all components isomorphic on $t$ edges. We therefore study whether this gap can be narrowed by considering the corresponding degree-associated edge-reconstruction numbers. Then we shall shift focus onto the degree-associated edge-reconstruction number of caterpillars and some other special classes of trees, with the main aim of obtaining some results towards determining the degree-associated edge-reconstruction number of a tree.


\section{Results on the edge-reconstruction number of a disconnected graph}

In \cite{m2}, Molina started to tackle the edge-reconstruction number of disconnected graphs. He showed that the edge-reconstruction results are similar to the vertex reconstruction results stated by Myrvold \cite {my1}, but a significant difference is that whereas the vertex reconstruction number of a graph is always three or more, the edge-reconstruction number of a disconnected graph is often two.
In summary, these are Molina's main results:

Let $G$ be a disconnected graph with at least four edges and at least two non-trivial components (that is, components that have more than one vertex). Then

\begin{itemize}
\item[(1)] if not all components are isomorphic, then $ern(G) \leq 3$;
\item[(2)] if all components are isomorphic, then $ern(G) \leq {t + 2}$ where $t$ is the number of edges in a component;
\item[(3)] if there exists a pair of non-isomorphic components in which one component has a cycle and $G$ does not have any components isomorphic to either $K_3$ or $K_{1,3}$, then $ern(G) \leq 2$.
\end{itemize}
He also observed that the value of $t + 2$ is attained, giving as an example the graph consisting of $p$ copies of $K_{1,t}$.

In \cite{al2}, Asciak and Lauri used line graphs in order to prove and extend Molina's results. In fact they proved the following results:
\begin{Theorem}
Let $G$ be a disconnected graph with at least four edges and the property that all components are isomorphic to a graph $H$. Then
\begin{itemize}
\item[(1)] if $H$ is isomorphic to $K_3$, then $ern(G) = 2$;
\item[(2)] if $H$ is isomorphic to $K_{1,3}$,  then $ern(G) = 5$;
\item[(3)] if $H$ is not isomorphic to $K_3$ or $K_{1,3}$, then $ern(G) \leq t + 2 $, where $t$ is the number of edges in $H$. Moreover, if $ ern(G) \geq {t + 1}$ then $H \simeq K_{1,t}$.
\end{itemize}
\end{Theorem}

\begin{Theorem} Let $G$ be a disconnected graph consisting of exactly two types of non-trivial components, namely those isomorphic to $K_3$ and those isomorphic to $K_{1,3}$. Then $ern(G) = 3$.
\end{Theorem}


They also tried to investigate conditions which force or do not allow $ern(G)$ to be equal to 2 and also showed that in general, there is no straightforward relationship between the edge-reconstruction number of $G$ and that of its components.

These results and data from Rivshin's computer search \cite{riv} (which showed that out of more than a billion graphs on at most eleven vertices, only fifty-six disconnected graphs have edge-reconstruction number greater than 3 and that out of these disconnected graphs, only four graphs do not have isolated vertices as  components, namely $2K_{1,2}, 2K_{1,3}, 2K_{1,4}$ and $3K_{1,2}$) led Asciak and Lauri to make the following conjecture.

\begin{Conjecture}
Suppose that $ern(G) > 3$ for a disconnected graph all of whose components are isomorphic to $H$. Then $H$ is isomorphic to the star $K_{1,r}$ where $r$ is the number of edges.
\end{Conjecture}

\section{Results on the edge-reconstruction number of a tree}

We shall first describe some special types of trees and also give some basic definitions on general trees.

A \emph{caterpillar} is a tree whose non-leaf vertices (a \emph{leaf} or an \emph{endvertex} is a vertex of degree 1) induce a path called the \emph{spine} of the caterpillar. A caterpillar will be represented as the sequence $\langle a_1, ..., a_n \rangle$ which denotes a caterpillar with spine $ v_1, ... , v_n $ such that $a_i$ leaf vertices are incident to $v_i$ for each $i \in \{1, ... , n\}$. Clearly this representation of a given caterpillar is unique up to left to right orientation. Note that for $ 2 \leq i \leq {n - 1} $, each $a_i \geq 0$, but $a_1, a_n \geq 1$. Such a sequence is called a \emph{caterpillar sequence}.

A special type of tree denoted by $S_{a,b,c}$ is a tree similar to a star (a star is the tree on $n$ vertices, $n - 1$ of which are endvertices) which consists of three paths on $a$, $b$ and $c$ edges, respectively, emerging from a common vertex. Some examples are shown in Figure 1.

\begin{figure}[h]
\centering
\includegraphics[scale=0.8]{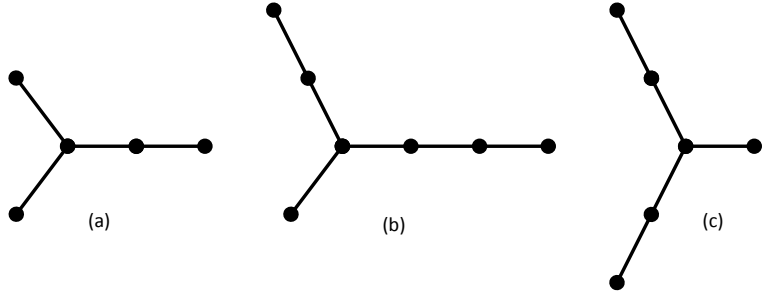}
\caption{The trees: (a) $S_{1,1,2}$; (b) $S_{1,2,3}$; (c) $S_{1,2,2}$}
\end{figure}

A tree is called a \emph{quasipath} if it is either a path $P_n$ on $n$ vertices or one of the two special trees $S_{1,1,2}$ and $S_{1,2,3}$.

We define the \emph{weight} of a vertex $v$ of a tree $T$, denoted by $wt(v)$, to be the number of vertices in a largest component of $T - v$. The \emph{centroid} of a tree $T$ is the set of all vertices with minimum weight; this weight denoted by $wt(T)$. A \emph{centroidal vertex} is a vertex in the centroid. It is well-known that the centroid of a tree consists of either one vertex or two adjacent vertices.  A tree with one centroidal vertex is called \emph{unicentroidal} while a tree with two centroidal vertices is called \emph{bicentroidal}. In the latter case, the edge joining the centroidal vertices is called the \emph{centroidal edge}.  When T is bicentroidal with centroidal edge $e$, then the two components of $T - e$ are said to be \emph{centroidal components}.

Rivshin \cite {riv}, using his computer program, obtained the following result which helps complete the theoretical results given in \cite {almp}.\\

\begin{Result}
\begin{itemize}
\item[(1)] If $ T$ is either one of the two trees $H_1$ and $H_2$ shown in Figure 2, then  $ern(T) = 3$.
\item[(2)] If  $T$ is a bicentroidal tree in which vertices $a$ and $b$ are the centroidal vertices so that the two centroidal components of $T - ab$ are $S_{1,2,3}$, then $ ern(T) = 2$. But if the two components are $S_{1,1,2}$ then  $ern(T) = 2$ unless $T$ is the tree $H_3$ shown in Figure 2 where $ern(H_3) = 3$.

\end{itemize}
\end{Result}
\begin{figure}[h]
\centering
\includegraphics[scale=0.8]{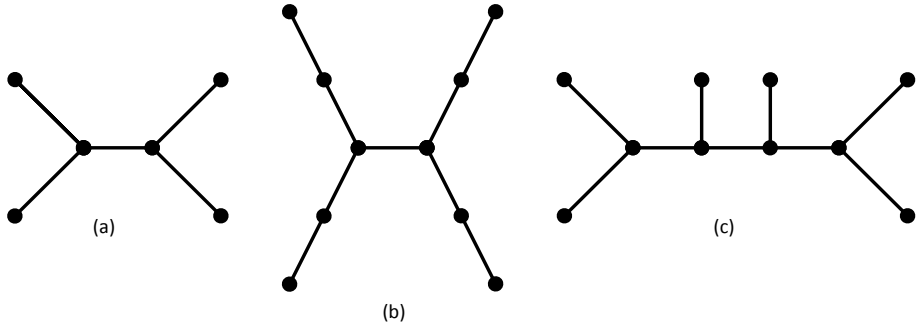}
\caption{The trees: (a) $H_1$ (b) $H_2$ (c) $H_3$}
\end{figure}
Molina \cite {m3} gave the following result.

\begin{Result} If $T$ is a unicentroidal tree on at least four edges then $\ern(T) \leq 3$
\end{Result}
Then Asciak, Lauri, Myrvold and Pannone \cite {almp} proved the following.
\begin{Result}  Every bicentroidal tree except the caterpillars with sequences $\langle 2,2 \rangle,\langle 2,1,1,2 \rangle$ and the graph $H_2$ shown in Figure 2 has edge-reconstruction number equal to 2.
\end{Result}


For unicentroidal trees, the results of the computer searches in \cite {almp} and those of Rivshin, led Asciak, Lauri, Myrvold and Pannone  to the infinite family of trees $T_k$ $(k \geq 2)$, where $k$ is the degree of the central vertex (depicted in Figure 3) having $ern = 3$ . Note that when $k = 2$, $T_k$ is the caterpillar $\langle 2,0,2 \rangle$.

\begin{figure}[H]
\centering
\includegraphics[scale=0.8]{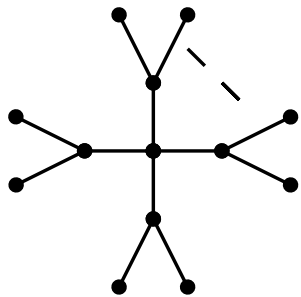}
\caption{An infinite family of trees $T_k$ with $ern = 3$}
\end{figure}

They also found the graph $G_{15}$ on fifteen vertices, shown in Figure 4, which does not fall within any known infinite class but which also has $ern = 3$. These computer searches and results presented in \cite {almp} led them to make the following conjecture for unicentroidal trees.

\begin{Conjecture}\cite {almp} The only infinite classes of trees which have $ern = 3$ are the paths on an odd number of vertices, the caterpillars $ \langle 2,0, ... ,0,2 \rangle $ of even diameter, and the family of trees $T_k$ described above.
\end{Conjecture}

\begin{figure}[h]
\centering
\includegraphics[scale=1]{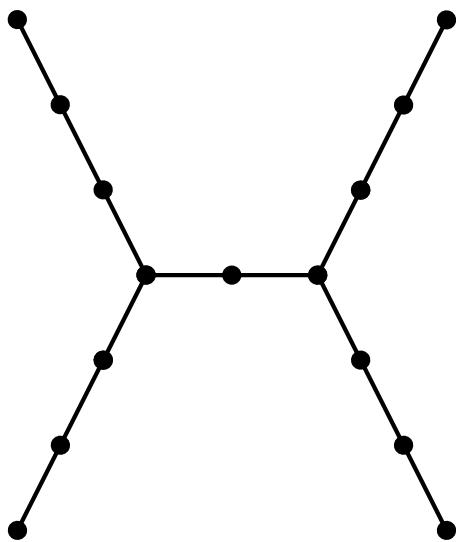}
\caption{The tree $G_{15}$ on fifteen vertices with $ern = 3$}
\end{figure}

\section{Degree-associated edge-reconstruction number of a disconnected graph}

Continuing on the work presented by Molina and more recently by Asciak and Lauri as described in Section 2, we shall shift our study to the degree-associated edge-reconstruction number of disconnected graphs wherein the degree of the deleted edge $d(e)$ is given together with the edge-card. But first we need the following two lemmas, the first of which is due to Ma et al. in \cite{mssw}.

\begin{Lemma}
If $G$ has an edge $e$ such that $d(e) = 0$ or no two non-adjacent vertices in $G - e$ other than the endpoints of $e$ have degree-sum $d(e)$, then da-ecard $(G - e, d(e))$ determines $G$.
\end{Lemma}

The condition in Lemma 4.1 is sufficient but not necessary for $(G - e, d(e))$ to determine $G$. In fact if $G$ is a graph in which an edge joins two disjoint complete graphs, then the condition fails, but $dern(G) = 1$.

\begin{Lemma}
Let $e$ be an edge in a graph $G$ and $(G - e, d(e))$ be a given da-ecard. If $G - e$ is without isolated vertices, then
\begin{itemize}
\item[(1)] if $d(e) = 2$, then $e$ is incident to two endvertices from $G - e$;
\item[(2)] if $d(e) > 2$, then $e$ is incident to at most one endvertex from $G - e$.
\end{itemize}
However, if the edge-card $G - e$ has isolated vertices and $d(e) \geq 2$, then the deleted edge can either join an isolated vertex to a non-endvertex or is incident to two endvertices whenever $d(e) = 2$ or  is incident to two non-endvertices whenever $d(e) > 2$.
\end{Lemma}

\noindent
{\bf Proof.}  Let $e $ be an edge in a graph $G$.

Suppose first that the edge-card $G - e$ has no isolated vertices. We show that (1) and (2) hold.
If $d(e) = 2$, then there must be two edges adjacent to $e$, so the only way to place the edge $e$ in order to obtain $(G - e, d(e))$ as a da-ecard is to join two endvertices.
If however $d (e) > 2$, then the missing edge $e$ can join either an endvertex to a non-endvertex or two non-endvertices.

Now suppose that the edge-card $G-e$ has at least one isolated vertex.
Then if $d(e) \geq 2$, there is the possibility (apart from the previous situation) that the missing edge joins an isolated vertex to a non-endvertex.\lesta \\

We are now in a position to find the degree-associated edge-reconstruction number of disconnected graphs for a number of cases which include those  mentioned in Section 2. But first we need the following definition.

\begin{Definition}  Suppose that a graph $H \not\simeq G$ has in its edge-deck the edge-cards $G - e_1, G - e_2, \ldots , G - e_k$, we then say that $H$ is a \emph{blocker} for these edge-cards or that $H$ \emph{blocks} these edge-cards.
\end{Definition}

\begin{Theorem}
Let $G$ be a disconnected graph with at least two non-trivial components all of which have at least three edges. Then
\begin{itemize}
\item[(1)] if all components are isomorphic to $K_3$, then $dern(G) = 1$;
\item[(2)] if all components are isomorphic to $K_{1,3}$,  then $dern(G) = 4$;
\item[(3)] if all components are isomorphic to  $K_{1,t}$, where $t$ is greater than 3, then $dern(G) = 1$;
\item[(4)] if all components are of exactly two types, namely those isomorphic to $K_3$ and those isomorphic to $K_{1,3}$, then $dern(G) = 2$;
\item[(5)] if $G$ has either

\begin{itemize}
\item[(i)] only $K_3$ components and at least one isolated vertex,
or

\item[(ii)] only $K_{1,3}$ components and at least one isolated vertex,

\end{itemize}
then in both cases $dern(G) = 4$.
\end{itemize}
\end{Theorem}

\noindent
{\bf Proof } \textbf{(1)}. As the graph $G$ is made up of only $p$ copies of the $K_3$ component, then its da-edeck consists of $3p$ copies of da-ecard $(G',2)$ where $G'$ has $(p - 1)$ copies of $K_3$ and a path $P_3$. By Lemma 4.1, every da-ecard ($G'$,2) determines $G$.

\smallskip\noindent
\textbf{(2)}. Since all components in graph $G$ are isomorphic to $K_{1,3}$ then all da-ecards of $G$ are copies of $(G^*,2)$ where $G^*$ consists of $(p - 1)$ copies of $K_{1,3}$, the path $P_3$ and an isolated vertex. By Lemma 4.2, the possible graphs having $(G^*,2)$ as a da-ecard is either the original graph $G$ or graph $H$ which is obtained from $G^*$ by joining the two endvertices of $P_3$. Since there exist only three da-ecards in the da-edeck of $G$ that are also in the da-edeck of $H$, therefore $dern(G) = 4$.

\smallskip\noindent
\textbf{(3)}. Since $G$ consists of $p$ copies of $K_{1,t}$ where $t$ is the number of edges, then all da-ecards are of the form $(G_1, t - 1)$ where $ G_1$ has $(p - 1) K_{1,t}$ components, a $K_{1,t - 1}$ component and an isolated vertex. By Lemma 4.1, every da-ecard $(G_1,t - 1)$ determines  $G$.

\smallskip\noindent
\textbf{(4)}. Let $G$ consist of $p$ copies of $K_{1,3}$ and $q$ copies of $K_3$, so the da-edeck of $G$ consists of two different da-ecards $(C_1,2)$ and $(C_2,2)$. The da-ecard $(C_1,2)$ is the graph $G - e_1$ where $e_1$ is an edge in $K_3$, so $C_1$ consists of $p$ copies of $K_{1,3}$, $(q - 1)$ copies of $K_3$ and a component $P_3$. The da-ecard $(C_2,2)$ is the graph $G - e_2$ where $e_2$ is an edge in $K_{1,3}$ and therefore $C_2$ consists of $q$ copies of $K_3$, $(p - 1)$ copies of $K_{1,3}$, a component $P_3$ and an isolated vertex. Suppose that graph $H_i \not\simeq  G$ can be obtained from $G - e_2$ by adding an edge $e$ whose $d(e) = 2$. By Lemma 4.2, $H_i$ must be one of the following graphs:

\begin{itemize}
\item $H_1$ is a disconnected graph whose components are a caterpillar $\langle 2,0,0,1 \rangle$, $(p - 2)$ copies of $K_{1,3}$, $q$ copies of $K_3$ and an isolated vertex,
\item $H_2$ is a disconnected graph whose components are a caterpillar $\langle 2,0,0,2 \rangle$, $(p - 3)$ copies of $K_{1,3}$, $q$ copies of $K_3$, a component $P_3$ and an isolated vertex,
\item $H_3$ is a disconnected graph whose components are a component $Z$, formed by adding an edge to two endvertices of a $K_{1,3}$ component,
$(p - 2)$ copies of $K_{1,3}$, $q$ copies of $K_3$, a component $P_3$  and an  isolated vertex.
\item $H_4$ is a disconnected graph whose components are component $Z$(a similar component to the above), which is obtained by joining an isolated vertex to a $K_3$ component, $(p - 1)$ copies of $K_{1,3}$, $(q - 1)$ copies of $K_3$ and a component $P_3$ .

\end{itemize}
In all cases there is at most one da-ecard in the da-edeck of $G$ that is in the da-edeck of $H_i$, $ 1 \leq i \leq 4 $. Thus $dern(G) = 2$ since graphs $H_1$, $H_2$, $H_3$ can only contain one da-ecard in common with $G$, while the specific choice $(C_1,2)$, $(C_2,2)$ of da-ecards rules out graph $H_4$.




\smallskip\noindent
\textbf{(5)}. Case (i): Let $G$ consist of at least two $K_3$ components and $r$ isolated vertices, $r \geq 1$. Then all da-ecards of $G$ are copies of $(G'',2)$ where $G''$ is a disconnected graph consisting of a $P_3$ component, $(p - 1)$ copies of $K_3$ and $r$ isolated vertices. A disconnected graph $F$ made up of $(p - 1)$ copies of $K_3$, a $K_{1,3}$ component and $(r - 1)$ isolated vertices has three copies of $(G'',2)$ in its da-edeck, so $dern(G) \geq 4$.

Using Lemma 4.2, the missing edge $e'$ can join the two endvertices in $G''$ in which case the resultant graph $H \simeq G $; otherwise, an isolated vertex is joined either to a $K_3$ component or $P_3$ component. In the latter cases, each resulting graph can share at most three da-ecards with $G$, so $dern(G) \leq 4$. Hence $ dern(G) = 4$.

\smallskip\noindent
Case (ii): Similar arguments presented in Case (i) hold for this case. \lesta

\smallskip\noindent
In the cases considered so far, we did not investigate the case when all the components of $G$ are all isomorphic and not equal to either $K_3$ or $K_{1,t}, t \geq 3$. The next results will address this situation.

\subsection{Disconnected graphs whose components are all isomorphic}

In this section we consider only $G = kH$, that is, the disjoint union of copies of the connected graph $H$, since otherwise as pointed out in Section 2, $ern(G)$ is at most 3 \cite{m2}. It is also known that if $t = |E(H)|$, $ern(G)$ can be as large as $t + 2$ and this can happen only when $H = K_{1,t}$ \cite{al2, m2}. However, in this case, $dern(G) = 1$ except for the case $H = K_{1,3}$ considered above. Therefore, can the degree-associated edge-reconstruction number help to reduce the gap between 3 and $t + 2$?

But first we need to define some definitions and state some important results which will be required throughout this sub-section.
%
%
A \emph{bipartite graph} is one whose vertex set can be partitioned into two sets such that no edge joins two vertices in the same set. These two sets are called the \emph{colour classes} of the bipartite graph, and they are uniquely defined if the graph is connected.

A graph $G$ is \emph{edge-transitive} if, given any two edges $\{a,b\}$ and  $\{c,d\}$,  there is an \emph{automorphism} $\alpha$ (that is, a one-one mapping of the vertex set of graph $G$ onto itself which preserves adjacency) such that $\{\alpha (a),\alpha (b)\} =  \{c,d\}$.
It is easy to show that if all edge-cards $H - e$ are isomorphic then $H$ must be edge-transitive.  The complete bipartite graph $K_{p,q}$ is an example of an edge-transitive graph. We shall also need the following easy result.

\begin{Proposition} Let $H$ be an $r$-regular graph.\\
Then $adv\text{-} dern(H) =  dern(H) = 1$.\lesta
\end{Proposition}


Recall that the minimum degree of a graph $G$ denoted by $\delta(G)$ is the smallest number of edges incident to any vertex $v$ in $G$.

We are now able to start considering the case when all the edge-cards of a disconnected graph are isomorphic.
\begin{Theorem}
Let $G$ be a disconnected graph, all of whose components are isomorphic to $H$, and suppose that all edge-cards of $H$ are isomorphic. Suppose also that $\delta(H) \geq 3$. Then $dern(G) \leq 2$.
\end{Theorem}

\noindent
{\bf Proof.} Since all edge-cards $H - e$ are isomorphic then $H$ is edge-transitive. Also all edges have degree $p = r + s - 2$, where $r \geq s$
are the degrees of the two endvertices of the edge; $r$ and $s$ are the same for all edges, that is, $H$ is either regular (when $r = s$) or bi-degreed (that is, all vertex degrees are either $r$ or $s < r$) although we do not know, from the da-ecards, the values of $r$ and $s$. Due to Proposition 4.1 we may assume that $H$ is not regular. In this case, $H$ is bipartite, with the vertices of different degree forming the two colour-classes of $H$.

We shall assume that $H$ has $a$ vertices of degree $r$ and  $b$ vertices of degree $s$, where $ar = bs$. Also, since we are assuming that $\delta(H) \geq 3$, then $s > 2$.

All da-ecards are of the type $( (H - e) \cup (t - 1)H, p)$.
Suppose we are given two such da-ecards $C_1, C_2$. Then we need to consider the following two cases.\\\\
\emph{Case 1}: $C_1, C_2$ are blocked by a graph $G' = (H - e) \cup K \cup (t - 2)H$, where $K \simeq H + f$ with $d(f) = p$ in $K$.

But since $d(f)$ must be $p$ in $K$ and its endvertices  can only have degrees $r$ and $s$ in $H$, the only possibility is that $f$ joins two vertices of degree $s$ in $H$, and therefore
\[ p = r + s - 2 = s + s \]
 \[ s = r - 2 \]

So, we have that $H $ has $a$ vertices of degree $r$ and $ b = \frac{ar}{s} $   vertices of degree $s$.  Also, $K$ has $a$ vertices of degree $r$, $b - 2$ vertices of degree $s$, and two vertices of degree $s + 1$.  Moreover, $K$ must have an edge $ f' \neq f$  such that $ K - f' \simeq H$   in order for $G'$ to be a blocker of the da-ecards $ C_1, C_2$. But $K - f'$ has a vertex of degree $s + 1 = r - 1$ which $H$ does not.

This contradiction shows that $ C_1, C_2$ cannot have a blocker of this type.\\\\
\emph{Case 2}: $C_1, C_2$ are blocked by a graph $G' = H' \cup (t - 1)H$, where $H' \not \simeq H $ but $H'$ has two da-ecards in common with $H$.

Therefore we can take $H'$ to be equal to $H - e + e'$ with $e = uv$, $d(u) = s$, $ d(v) = r$ in $H$ and $e' \neq e$, but such that $d(e')$ in $H'$ is equal to  $d(e) = r + s - 2$ in $H$.  But $d(e')$ in $H'$ is equal to $ r_1 + r_2$ where $r_1, r_2 \in \{r, r - 1, s, s - 1\}$ (these being the degrees of the vertices in $H - e$) but  $\{r_1, r_2 \} \neq \{r - 1, s - 1\}$ (since $ e' \neq e$).

So we have these two possibilities:

\begin{itemize}
\item[(a)] $e'$ is adjacent to two vertices of degree $s$ in $H - e$. Therefore
\begin{eqnarray*}
 r + s - 2 &=& s + s\\
  s &=& r - 2
\end{eqnarray*}

\item[(b)] $e'$ is incident to two vertices of degree $s$ and $s - 1$ in $H - e$. Therefore
\begin{eqnarray*}
 r + s - 2 &=& s + s - 1\\
  s &=& r - 1
\end{eqnarray*}
\end{itemize}

We shall now consider each of these two possibilities:\\
\emph{Case 2(a)}:
Let $e' = xy$. Therefore $H'$ has one vertex $u$ of degree $s - 1$, the vertex $v$ of degree $r - 1$ and the vertices $x, y$ of degree $s + 1$. Since $s= r - 2$ then $r - 1 = s + 1$, and $x, y, v$ are the only vertices of $H'$ with this degree. Also, since $s \geq 3$, each of the two vertices $x, y$ are adjacent, in $H$, to at least three vertices of degree $r$. Therefore, in $H'$, each of $x, y$ is adjacent to at least two vertices of degree $r$.
But we know that for $G'$ to block both $C_1$ and $C_2$, $H'$ must have an edge $f \neq e'$ such that $H' - f$ is isomorphic to an edge-card of $H$.   Comparing the degrees in an edge-card of $H$ with the degrees in $H'$ we see that $f$ must be adjacent, in $H'$, to two vertices of degree $s + 1$.  But since $f \neq e$, the endvertices of $f$ must be $v$ and one of $x$ or $y$. This means that in $H' - f$ there still is a vertex ($x$ or $y$) of degree $s + 1$ adjacent to a vertex of degree $r$. But $H - e$ contains no such pair of adjacent vertices.  Therefore, this case cannot happen.\\\\
\emph{Case 2(b)}:
We proceed similarly to the previous case. Let $e' = xy$ and suppose that $d(x) = s - 1$ and $d(y) = s$ in $H - e$.  Therefore $x$ must be the vertex $u$ and $y$ cannot be the vertex $v$, since $e \neq e'$.  Therefore, in $H'$, $y$ is a vertex of degree $s + 1 = r$  adjacent to one vertex, $x$, of degree $s$, possibly another vertex, $v$, also of degree $s$, and the remaining $r - 2$ neighbours all of degree $r$. Since $ s \geq 3$ and $r = s + 1$, this means that $y$ is adjacent to at least two non-adjacent vertices of degree $r$.

As before, for $G'$ to be a blocker of $C_1$ and $C_2$, $H'$ must have an edge $f \neq e'$  such that $(H' - f, d(f))$ is the same as a da-ecard of $H$. Therefore $f$ must be incident to a vertex of degree $s$ and a vertex of degree $r$.  Now, since $f \neq e'$, these two vertices cannot be $x$ and $y$.  Also, $f$ cannot be $yv$ (supposing $v$ were adjacent to $y$ in $H$), because in this case, $v$  would be a vertex of degree $s-1$ adjacent to vertices of degree $s$ in $H' - f$, which no edge-card of $H$ has. But then, since $y$ is adjacent to at least two non-adjacent vertices of degree $r$ in $H'$, then it must be adjacent to at least one vertex of degree $r$ in $H' - f$, that is, $H' - f$ contains an edge joining two vertices of degree $r$, which no edge-card of $H$ has. Therefore even here we have shown that this case cannot hold.\lesta\\\\
\textbf{Remark.}  If we allow $\delta(H) = 1$ and $ k > 1$, then $G = kK_{1,3}$ has $dern$ equal to 4 and $dern(kK_{1,2}) = 3$ (see later) and if we allow $\delta(H) = 2$ then $G = kK_{2,3}$ has $dern$ equal to 3. In fact we believe that these are the only three cases when $dern(G) > 2$, and we conjecture:
\vskip30pt

\begin{Conjecture}
Let $G$ be as in Theorem 4.2 and suppose that the condition $\delta(H) \geq 3$  is replaced by $H \neq K_{1,3}$, $H \neq K_{1,2}$ and $H \neq K_{2,3}$. Then $dern(G) \leq 2$.
\end{Conjecture}

\smallskip\noindent
We now consider the case when not all the edge-cards of $H$ are isomorphic. It seems that this case is more difficult to settle. So we shall give only  a very partial result.

But first we need to define the minimum multiplicity of a graph. The \emph{multiplicity} of an edge-card of $H$ is the number of times it appears in the edge-deck of $H$, and the \emph{minimum multiplicity} of $H$, denoted by $mm(H)$ is the  minimum amongst all multiplicities of edge-cards appearing in the edge-deck of $H$.

We give this simple result which does not use the concept of $dern(H)$. We believe that the result is far from being the best possible.

\begin{Theorem} Let $G = kH$ and suppose not all edge-cards of $G$ are isomorphic. Then $ ern(G) \leq min\{adv\text{-} ern(H), 2 + mm(H)\}$
\end{Theorem}
\noindent
{\bf Proof.} Let $H_1, H_2$ be two non-isomorphic edge-cards of $H$. Then, any collection of edge-cards of $G$ which contains $ H_1 \cup (k - 1)H$ and  $ H_2 \cup (k - 1)H$ can only be blocked by a graph $G' = H' \cup (k - 1)H$. In this case, if the collection of edge-cards is  $\{ H_1 \cup (k - 1)H, \ldots , H_r \cup (k - 1)H\}$ then $H'$ must block the edge-cards $H_1, H_2, \ldots , H_r$.

If these edge-cards are all obtained from the same component of $G$, and $r \geq adv\text{-} ern(H) $ then no $H'$ can block them.

If we let $H_2, \ldots , H_r$ be isomorphic copies of a single edge-card of $H$ and $ r = 1 + mm(H)$ then again no $H'$ can block $H_1, H_2, \ldots , H_r$. \lesta

In the following section we shall investigate $dern(G)$, where $G$ is a special case of disconnected graphs whose components are all isomorphic but not all of its edge-cards are isomorphic.

\subsection{Disconnected graphs whose components are isomorphic to $P_n$}

We can easily verify that some special trees such as paths $P_n$ and stars $K_{1,n}$ have degree-associated edge-reconstruction number of 1.

In Theorem 4.1 we have shown that even $k$ copies of stars $K_{1,n}$ for $n > 3 $ have degree-associated edge-reconstruction number equal to 1. But we shall show that while $dern(kP_n) = 2 $ for $ n > 3$,  $dern(kP_3) = 3$.

\begin{Example}
 $dern(kP_n) = 2$ for $ n > 3$ and $k > 1$
\end{Example}
\noindent
Let $G$ be a graph consisting of $k$ copies of $P_n$. Then there are two types of da-ecards which are either $C_1 = (S , 1)$, where $S$ is a disconnected graph having $k - 1$ components isomorphic to $P_n$, a $P_{n - 1}$ component and an isolated vertex, or $C_2 = (R, 2)$, where $R$ is a disconnected graph having $k - 1$ components isomorphic to $P_n$  and two other path components $P_l$ and $P_{n-l}$ where $l < n$. Now graphs other than $G$ having either da-ecard $C_1$ or $C_2$ are $G_1 = P_{n - 1} \cup P_{n + 1} \cup (k - 2)P_n$ and $G_2 = P_{n - l} \cup P_{n + l} \cup (k - 2)P_n$ where $l < n$, respectively. This means that $dern(kP_n) > 1$, so $C_1$ or $C_2$ alone do not give $kP_n$ uniquely where $n > 3$.

We now show that $C_1$, $C_2$ determine $G$.
\smallskip\noindent
Let $G'$ be a reconstruction from \{$C_1,C_2$\} and suppose that $G'$ is obtained by adding a new edge $e'$ to the edge-card $S$. Now since the degree of $e'$ is 1, the isolated vertex can join an endvertex of a component $P_n$ to obtain the graph $G_1$. But the resultant graph $G_1$ does not have the edge-card $R$. So the only possibility is that $G'$ is obtained by joining the isolated vertex to an endvertex of path $P_{n-1}$, so $G'\simeq G$. \lesta

\begin{Example}
 $dern(kP_3) = 3$
\end{Example}
\noindent
Let $G$ be a graph consisting of $k$ copies of $P_3$. Then $G$ has only one distinct da-ecard which is of the form $(S,1)$ where $S$ is a disconnected graph having $k - 1$ copies of component $P_3$, and components $K_2$ and $K_1$. With just one such da-ecard one can conclude  that any reconstruction has paths for all its components. Now the graph $G_1 = P_2 \cup P_4 \cup (k-2)P_3$ has two da-ecards in common with $G$ and so two da-ecards are not sufficient to reconstruct $G$ uniquely. But three da-ecards  suffice in order to force a reconstruction to have all components isomorphic to $P_3$ because the reconstruction can only be obtained from $S$ by joining the isolated vertex to the smallest component $P_2$.\lesta

\section{Degree-associated edge-reconstruction number of a tree}


The empirical evidence provided by David Rivshin showed that, after investigating graphs on at most eleven vertices, only seventeen trees have edge-reconstruction number equal to 3. Three out of these seventeen trees are the bicentroidal trees $H_1$, $H_2$ and $H_3$ shown in Figure $2$ . We can easily check directly that $dern (H_2) = dern (H_3) = 1 $ while $dern (H_1) = 2 $. The remaining trees are unicentroidal and we can also show by hand that, the caterpillars $\langle 2,0,2 \rangle$,   $\langle 2,1,2 \rangle$, $\langle 2,3,2 \rangle$, the paths of odd order $P_5$, $P_7$,$P_9$ and $P_{11}$ and the tree $G_1$ shown in Figure 3 have $dern  = 1 $ while the remaining six unicentroidal trees, namely,  $\langle 1,0,1,0,1)$, $\langle 2,0{^3},2)$,  $\langle 1,0,1,0,1,0,1  \rangle$, $\langle 2,0{^5},2 \rangle$ , $S_{2,2,2}$ and $S_{3,3,3}$ can also be directly checked by hand to show that their $dern = 2 $. Previously we have stated that from computer search, the edge-reconstruction number of graph $G_{15}$ of Figure 4 is 3. But also, in this case, we can show directly that $dern(G_{15}) = 2$.

\subsection{Caterpillars}

Barrus and West \cite{bw} have shown that, except  for the case of the 6-vertex caterpillar $ H_1$ shown in Figure 2, $drn$ of a caterpillar is 2, or 1, for stars.

We shall now restrict our study in this section to the degree-associated edge-reconstruction number of caterpillars. But we first quote a very useful result by Molina \cite {m2} which allows us to identify a graph as a tree from two given edge-cards.
\begin{Lemma} Let $G$ be a graph with edges $e_1$ and $e_2$. Suppose that edge-card $G - e_1$ has two components which are trees of orders $a_1$ and $a_2$ while edge-card $G - e_2$ has another two components which are trees of orders $b_1$ and $b_2$. If $\{ a_1, a_2\} \neq \{b_1, b_2\}$, then $G$ is a tree.
\end{Lemma}

Moreover, we shall also make use of the next observation in order to identify when a tree is a caterpillar.

\begin{Observation} Let $T$ be a tree and let $e$ be an edge whose degree $d(e)$ is greater than 1 and that $T - e$ is a caterpillar with an isolated vertex. Then $T$ is a caterpillar.
Therefore the fact that $T$ is a caterpillar, is recognisable from any two da-ecards as long as at least one of them is not obtained by deleting an edge which changes the spine of the caterpillar.
\end{Observation}


\noindent
From now on, we shall assume that $T$ is a caterpillar and since a caterpillar sequence $\langle 1,0, \ldots , 0,1 \rangle$ denotes a path which has $dern = 1$ we shall henceforth assume that $T$ is not a path.
We shall only consider da-ecards corresponding to $T - e$ and $ T - f$ where $ e$ and $f$ are end-edges which do not change the spine of $T$.
This is equivalent  to reconstructing a caterpillar sequence $\langle a_1, \ldots , a_n\rangle$ from two sequences $\langle a_1, \ldots , a_i - 1, \ldots , a_n \rangle$ and $\langle a_1, \ldots , a_j - 1, \ldots , a_n \rangle$ and we shall consider the problem this way. Note that the two sequences are given only up to left to right orientation; we call two such orientations the reverse of each other. We shall also call  $\langle a_1, \ldots , a_i - 1, \ldots , a_n \rangle$ a \emph{reduction} of the sequence  $\langle a_1, \ldots , a_n\rangle$. A reduction is therefore a description of the corresponding da-ecard of $T$ without the isolated vertex. If a number of reductions of a sequence determines the sequence, we can say that the reductions reconstruct the sequence.  We shall need the definition of conjugate entries in a caterpillar sequence.
Given the caterpillar sequence $\langle a_1,a_2, \ldots , a_n \rangle$, then the entries $a_i$, $a_{n - i + 1}$ are said to be \emph{conjugate}.\\\

\noindent
{\bf Notation.} The term $a_i^-$ will stand for $a_i - 1 $ and  $a_i^+$ for $a_i + 1 $.

\begin{Lemma} If the caterpillar sequence $S = \langle a_1,a_2, \ldots , a_n \rangle$ is not reconstructible from the reductions $\langle a_1, \ldots , a_i^- , \ldots , a_n \rangle$ and $\langle a_1,\ldots , a_{n - i + 1}^- , \ldots , a_n \rangle$ then $ a_i = a_{n - i + 1}$ and all other conjugate entries are equal, except for a pair which differ exactly by one.
\end{Lemma}
\noindent
{\bf Proof.} For non-reconstructibility we must be able to get another sequence from the alignment

\[\langle a_1, \ldots , a_i^-, \ldots ,   a_{n - i + 1} , \ldots , a_n \rangle\]
and
\[\langle a_n, \ldots ,  a_{n - i + 1}^-  , \ldots  a_i , \ldots , a_1 \rangle.\]
Suppose, for contradiction, that $a_i \not= a_{n - i + 1}$ , therefore $ a_i^- \not = a_{n - i + 1}^-$.
Therefore for the above alignment to lead to a reconstruction of $S$, all other conjugate entries must be equal and $a_i , a_{n - i + 1}$ must differ by exactly one, say, without loss of generality, $a_i = a_{n - i + 1}^+$.
Then the given alignment reconstructs as
\[\langle a_n, a_{n - 1}, \ldots ,  a_{n - i + 1}, \ldots , a_i,  \ldots , a_2, a_1 \rangle.\]
But, since all conjugate entries apart from $a_i,  a_{n - i + 1}$ are equal, this sequence is simply the reverse of the original sequence $S$, therefore reconstruction is unique.
(The above reasoning can be noted by following example $5.1$ ).

Hence, for non-unique reconstruction, $a_i =  a_{n - i + 1}$ and $a_i^- = a_{n - i + 1}^-$. Therefore, for the above alignment to lead to a reconstruction of $S$, we must have for some $j$, that the conjugate pair $a_j,  a_{n - j + 1}$ differ by exactly one, say,  $a_j = a_{n - j + 1}^+$, and all other conjugate entries are equal, as required. Thus we get the following reconstruction
\[\langle a_1, \ldots , a_i^-, \ldots , a_j, \ldots ,  a_{n - j + 1}^+, \ldots , a_{n - i + 1}, \ldots ,a_n \rangle .\]
(Example $5.2 $ illustrates this situation). \lesta\\

\begin{Theorem} A caterpillar sequence can be reconstructed from two reductions. Therefore, if $T$ is a caterpillar with at least two end-edges whose removal does not change its spine, then $dern(T) \leq 2$.
\end{Theorem}
\noindent
{\bf Proof.} Suppose the caterpillar sequence $\langle a_1, \ldots , a_n \rangle $ is not reconstructed from $\langle a_1^-,\ldots , a_n \rangle $ and $\langle a_1, \ldots , a_n^- \rangle$. Then by Lemma 5.2, $a_1 = a_n$ and, for some $j$, the conjugates $a_j, a_{n - j + 1}$ differ by one. But then, since  $a_j, a_{n - j + 1}$ are not equal, the sequence is reconstructed from $ \langle a_1, \ldots , a_j^- , \ldots , a_n \rangle$ and $ \langle a_1, \dots , a_ {n - j + 1}^- , \ldots , a_n \rangle$.\lesta

\begin{Corollary} If $T$ is a caterpillar which is not a path, then $dern(T) \leq 2$.
\end{Corollary}
\noindent
{\bf Proof.} The only remaining case to consider is when the caterpillar $T$ has only one end-edge whose removal does not change its spine. But it is easy to check that for such caterpillars $dern(T) \leq 2$.\lesta

\begin{Example} A caterpillar is expressed by the sequence $ \langle 3,4,3,7,7,2,4,3 \rangle $ and the following two reductions representing two edge-deleted caterpillars  are
\[\langle 3,4,2,7,7,2,4,3 \rangle\]
\[\langle 3,4,1,7,7,3,4,3 \rangle\]
where the second reduction is reversed with respect to the first one.
Comparing the two reductions will give the sequence $\langle 3,4,2,7,7,3,4,3 \rangle$, which is simply the reverse of the original sequence of the caterpillar.
\end{Example}
\begin{Example}
Let a caterpillar be expressed by the sequence $\langle 2,7,3,5,3,6,2 \rangle$ and the deletion of two of its end-edges gives the following reductions
\[\langle 1,7,3,5,3,6,2 \rangle\]
\[\langle 1,6,3,5,3,7,2 \rangle\]
in which the second reduction is reversed.
By comparing the two reductions, the  sequence $\langle 1,7,3,5,3,7,2 \rangle$ is obtained, which is an alternative sequence to the original one. So the two sequences given by the two corresponding edge-cards do not reconstruct the given caterpillar.
\end{Example}

\subsection{Star-like trees}

In \cite {mdr}, Monikandan et al. defined a \emph{star-like tree} to be a tree obtained from $K_{1,m}, m \geq 3$ by subdividing each edge at most once. They proved that any such tree has  $dern$ at most 2.

As already pointed out in Section 4.1, stars $K_{1,n}$ have degree-associated edge-reconstruction number of 1.
We shall show that if $S_{p+1}^n$ where $n > 1, p > 0$ is a tree obtained from the star $K_{1,n}$  by subdividing each edge $p$ times, as shown in Figure 5, then  $dern(S_{p+1}^n) = 2 $.

\begin{figure}[h]
\centering
\includegraphics[scale=1]{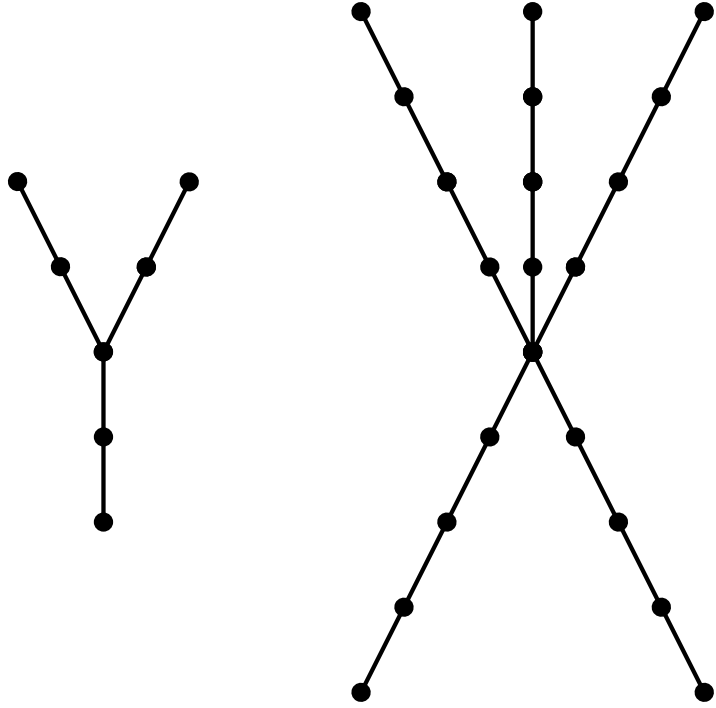}
\caption{The trees: (a)$S_{2}^3$  (b) $S_{4}^5$}
\end{figure}

\vskip20pt

\begin{Example}
Let $S_{p+1}^n$ be a tree defined as above. Then $dern(S_{p+1}^n) = 2$.
\end{Example}
\noindent
By considering all possible da-ecards of $S_{p+1}^n$, it is easy to show that it is not possible to reconstruct $S_{p+1}^n$ uniquely from only one da-ecard, therefore $dern(S_{p+1}^n) > 1$.
Let two da-ecards of the form $(C^*, 1)$ each be obtained by deleting an end-edge from $S_{p+1}^n$, so $C^*$ is made up of a unicentroidal tree in which all paths that emerge from the centroidal vertex are of the same order except one path whose order is one less than the others, and an isolated vertex.

We claim that these two da-ecards reconstruct $S_{p+1}^n$ uniquely.
Since in both cases the degree of the deleted edge is 1, the missing edge can only join the isolated vertex to any of the endvertices of $C^*$. But if the isolated vertex is joined to any one of the equal paths that emerge from the centroidal vertex then the resultant tree cannot have the two above mentioned da-ecards as part of its degree-associated edge-deck. Hence the only possible alternative is to join the isolated vertex to the smallest path emerging from the centroidal vertex in $C^*$. But this is isomorphic to $S_{p+1}^n$.
Hence $dern(S_{p+1}^n) = 2$.\lesta

Based on these findings, the results in Section 3, and the fact that  $dern(G) \leq ern(G)$, we offer the following conjecture.
\begin{Conjecture} If $T$ is a tree, then $dern(T) \leq 2 $.
\end{Conjecture}

\section*{Acknowledgement}
I am grateful to Josef Lauri for his many helpful suggestions which largely improved the style of this paper and Sections 4.1 and 5.

\end{document}